\documentclass{amsart}
\usepackage{tipa}
\usepackage{amssymb}
\usepackage{stmaryrd}
\usepackage{amsmath}
\usepackage{graphicx}

\newtheorem{theorem}{Theorem}[section]

\newtheorem{conjecture}[theorem]{Conjecture}

\theoremstyle{definition}

\theoremstyle{remark}

\numberwithin{equation}{section}

\begin{document}

\title{Figurate primes and Hilbert's 8th problem}
\author{Tianxin Cai}
\address{Department of Mathematics, Zhejiang University,
Hangzhou, 310027, People's Republic of China }
\email{txcai@zju.edu.cn}
\thanks{Project supported by the Natural Science Foundation of China~(No.~11351002).}


\author{Yong Zhang}
\address{Department of Mathematics, Zhejiang University,
Hangzhou, 310027, People's Republic of China }

\email{zhangyongzju@163.com}

\author{Zhongyan Shen}
\address{Department of Mathematics, Zhejiang International Study University,
Hangzhou, 310012, People's Republic of China }

\email{huanchenszyan@yahoo.com.cn}


\date{}


\begin{abstract}
In this paper, by using the theory of elliptic curves, we discuss
several Diophantine equations related with the so-called figurate
primes. Meanwhile, we raise several conjectures related with
figurate primes and Hilbert's 8th problem, including Goldbach's
conjecture, twin primes conjecture and Catalan's conjecture as well.

\end{abstract}

\maketitle
\section{Figurate primes}

In a letter to Crelle's Journal in 1844, Catalan stated that 8 and 9
are the only consecutive perfect powers, i.e., the Diophantine
equation
\[p^a-q^b=1\]
has unique positive integral solution $(p,q,a,b)=(2,3,3,2).$  This
is later known as Catalan's conjecture.

In 2004, P. Mih$\check{a}$ilescu \cite{Mih} proved this conjecture
by making extensive use of the theory cyclotomic fields and Galois
modules.

More generally, we have the Diophantine equation
\begin{equation}
p^a-q^b=k,
\end{equation}
where $p,q$ are primes, $a,b,k\in\mathbb{Z}^+.$

When $a=1,q=2,k=1,$ the solutions of $(1.1)$ are exactly Fermat primes.

When $p=2,b=1,k=1$, the solutions of $(1.1)$ are exactly Mesernne primes.

When $a,b\geq 2,k=1,$ (1.1) is the Diophantine equation for Catalan's conjecture.

For $k>1$, there are many authors who have investigated this
problem, more information can be found in \cite{Guy}: D9 Catalan
conjecture \& Difference of two powers and D10 Exponential
diophantine equations.

When $a=b=1,k=2$, (1.1) is the Diophantine equation for twin primes conjecture.

In spring 2013, the first author \cite{Cai} defined figurate primes
as the positive binomial coefficients
\[\binom{p^a}{i},~a\geq 1,~i\geq 1,\]
where $p$ is a prime. The set includes all primes, but with the same
density as the set of primes. We study the following Diophantine
equation
\begin{equation}
\binom{p^a}{i}-\binom{q^b}{j}=k,
\end{equation}
where $p,q$ are primes, $a,b,i,j,k\in \mathbb{Z}^+.$

When $k=1$, for $j=1,i\geq2$, we use elementary method to prove

\begin{theorem}\label{thm11}
For $(i,j)=(2,1)$, $(1.2)$ has exactly four solutions $(p, q, a, b)
= (2, 5, 2, 1),(3, 2,1, 1),(2, 3, 3, 3),(5, 3, 1, 2).$ For
$(i,j)=(3,1)$, $(1.2)$ has exactly three solutions $(p, q, a, b) =
(2, 3, 2, 1),(3, 83, 2, 1),(5, 3, 1, 2).$ For $(i,j)=(4,1)$, $(1.2)$
has exactly two solutions $(p, q, a, b) = (5, 2, 1, 2),(3, 5, 2,
3).$ \end{theorem}

For $i=b=1, j= 2.$ If $a$ is even, it's easy to see that (1.2) has
unique solution $(p,q)=(2,3);$ if $a=1,$ it seems likely that (1.2)
has infinitely many solutions, i.e., there are infinitely many pair
of primes $(p,q)$ satisfying \[p-1=\binom{q}{2}.\]
However, it is even a harder problem than that of prime representations by binary forms. The least 10 examples are $(p,q)=(2,2),(11, 5),(79,13),(137,17),(821,41),$\\
$(1831,61),(3917, 89),(4657, 97),(5051,101),(6329,113);$ if $a>1$ is
odd, we guess that (1.2) has no solutions. It's true for $a=3$ by an
easy calculation.

Similarly, as the proof of Theorem 1, we can get that all the solutions of the Diophantine equation
\[p^a-1=\binom{q^b}{3}\]
are $(p,q,a,b)=(5,2,1,2),(2,3,1,1),(11,5,1,1).$

For $i=j\geq2$, it's easy to verify that $(1.2)$ has no solutions.
By using the theory of elliptic curves, we have

\begin{theorem}\label{thm11}
For $(i,j)=(2,3)$, $(1.2)$ has unique solution $(p,q,a,b)=(3,7,2,1)$; and for $(i,j)=(3,2)$, $(1.2)$ has exactly two solutions $(p,q,a,b)=(2,3,2,1),(3,7,2,$\\
$1)$. For $(i,j)=(2,4)$, $(1.2)$ has exactly two solutions $(p,q,a,b)=(2,5,2,1),(3,7,2,$\\
$1)$; and for $(i,j)=(4,2)$, $(1.2)$ has no solutions.
\end{theorem}

When $k=2$, we have

\begin{theorem}\label{thm11}
For $(i,j)=(2,3)$, $(1.2)$ has exactly two solutions
$(p,q,a,b)=(2,2,2,2),(3,3,1,1)$; and for $(i,j)=(3,2)$, $(1.2)$ has
no solutions. For $(i,j)=(2,4)$, $(1.2)$ has unique solutions
$(p,q,a,b)=(3,2,1,2)$; and  for $(i,j)=(4,2)$, $(1.2)$ has unique
solution $(p,q,a,b)=(5,3,1,1)$.
\end{theorem}

\vskip20pt
\section{Hilbert's 8th problem}
Among the 23 problems that David Hilbert raised at the International
Congress of Mathematicians in Paris in 1900, the 8th one might be
the most profound and difficult, it includes Riemann Hypothesis,
Goldbach's conjecture and twin primes conjecture.

Goldbach's conjecture (1742)  is one of the most important unsolved
problems in number theory:

Every even integer greater than 2 can be expressed as the sum of two
primes, i.e.,
\[n=p+q,~n\geq4,\]
where $n$ is even and $p,q$ are primes. And every odd integer
greater than 5 can be expressed as the sum of three primes, i.e.,
\[n=p+q+r,~n\geq7,\]
where $n$ is odd and $p,q,r$ are primes. The first half for even
number is still an open problem.

However, by fundamental theorem of arithmetic, each positive integer
can be constructed from the product of primes, prime numbers are the
basic building blocks of any positive integer in multiplication.
Meanwhile, it seems that primes don't play key role in addition.
Besides, it's not a perfect result that each even integer is the sum
of two primes while each odd integer is the sum of three primes,
according to Goldbach's conjecture.

What we have to point out is that, among the composites in figurate
primes the amounts of even integers is as many as those of odd
integers, they are many more than powers of 2, cf. \cite{Guy}: A19
Values of $n$ making $n-2^k$ prime \& Odd numbers not of the form
$\pm p^a\pm2^b$. By calculations with computer, we check and find
that every positive integer $1<n\leq10^7$ can be expressed as the
sum of two figurate primes, i.e., the Diophantine equation
\begin{equation}
n=\binom{p^a}{i}+\binom{q^b}{j}
\end{equation}
always has solutions with primes $p,q$ and $a,b,i,j \in
\mathbb{Z}^+$. Therefore, we raise the following

\begin{conjecture}\label{conj21}
Every positive integer $n>1$ can be expressed as the sum of two
figurate primes.
\end{conjecture}

\begin{conjecture}\label{conj22}
$($week twin primes conjecture$)$ There are infinitely many pairs of
figurate primes with difference $2$.
\end{conjecture}

Moreover, we read Hilbert's speech carefully, in the satement of the
8th problem he mentioned

\begin{quote}\emph{``After an exhaustive discussion of Riemann's prime number
formula, perhaps we may sometime be in a position to attempt the
rigorous solution of Goldbach's problem, viz., whether every integer
is expressible as the sum of two positive prime numbers; and further
to attack the well-known question, whether there are an infinite
number of pairs of prime numbers with the difference 2, or even the
more general problem, whether the linear diophantine equation
                \[ax + by + c = 0\]
(with given integral coefficients each prime to the others) is
always solvable in prime numbers $x$ and $y$."}\end{quote}

With the idea of figurate primes and by numerical calculations, we
have

\begin{conjecture}\label{conj23}
For any positive integers $a$ and $b$, $(a, b)=1,$ when~$n\geq
(a-1)(b-1)$, there always exists prime pair $(x, y)$, such that
\[ax + by = n.\]
\end{conjecture}

Meanwhile, if we call a positive integer $n$ a proper one if $n$ is
a figurate prime but not a prime. Then we even have a stronger

\begin{conjecture}\label{conj24}
Every positive integer $n>5$ can be expressed as the sum of a prime
and a proper figurate prime.
\end{conjecture}

\vskip20pt
\section{Proof of the Theorems}
{\it \textbf{Proof of Theorem 1.1.}} When $k=1$, for $(i,j)=(2,1),$
$(1.2)$ is equal to \[(p^{a}+1)(p^{a}-2)=2q^{b}.\] Clearly,
$d=(p^{a}+1,p^{a}-2)=1$ or $3$.

If $d=1,p=2$, we have
\begin{equation*}
2^{a}-2=2,2^{a}+1=q^{b}
\end{equation*}
or
\begin{equation*}
2^{a}-2=2q^{b},2^{a}+1=1,
\end{equation*}
it's easy to see that $(p, q, a, b) = (2, 5, 2, 1)$.

If $d=1,p>2$, we have
\begin{equation*}
p^{a}-2=q^{b},p^{a}+1=2
\end{equation*}
or
\begin{equation*}
p^{a}-2=1,p^{a}+1=2q^{b},
\end{equation*}
 it's easy to see that $(p, q, a, b) = (3, 2, 1, 1)$.

If $d=3$, then $q=3$. For $p=2$, we have
\begin{equation*}
2^{a}-2=6,2^{a}+1=3^{b-1}
\end{equation*}
or
\begin{equation*}
2^{a}-2=2\cdot3^{b-1},2^{a}+1=3,
\end{equation*}
it's easy to see that $(p, q, a, b) = (2, 3, 3, 3)$.

For $p>2$, we have
\begin{equation*}
p^{a}-2=3,p^{a}+1=2\cdot3^{b-1}
\end{equation*}
or
\begin{equation*}
p^{a}-2=3^{b-1},p^{a}+1=6,
\end{equation*}
it's easy to see that $(p, q, a, b) = (5, 3, 1, 2)$.

As for the cases $(i,j)=(3,1)$ or $(4,1)$, we can prove by using similar method. \hfill $\Box$\\

{\it \textbf{Proof of Theorem 1.2.}}  When $k=1$, for convenience,
put $p^a=y,q^b=x$ for $(i,j)=(2,3)$ and $p^a=x,q^b=y$ for
$(i,j)=(3,2)$ in (1.2), respectively. Let
\[x=\frac{X+12}{12},~y=\frac{Y+36}{72},\] the converse
transformation is
\[X=12x-12,~Y=36(2y-1).\]
Then, we have
\[Y^2=X^3-144X+11664,~Y^2=X^3-144X-3024,\]respectively.

Using \textbf{Magma}, we get all the integral points on the above two elliptic curves. For $Y^2=X^3-144X+11664$, the point $(X,Y)=(72,612)$ leads to the unique solution $(p,q,a,b)=(3,7,2,1)$ of (1.2). For $Y^2=X^3-144X-3024$, the points $(X,Y)=(36,180),(84,756)$ lead to the two solutions $(p,q,a,b)=(2,3,2,1),(3,7,2,1)$ of (1.2).

Let $p^a=y,q^b=x$ for $(i,j)=(2,4)$ in (1.2), let
\[x=\frac{X+3}{6},~y=Y+2,\] we have
\[Y^2=3X^4+6X^3-3X^2-6X+81.\]

Using \textbf{Magma}, we get all the integral points on this elliptic curve, they are
\[(X,\pm Y)=(10,-11;189),(-1,-2,1,0;9),(-4,3;-21),(-6,5;51),(-92,91;14499),\]
by some calculations, we find that $(X,Y)=(3,21),(5,51)$ lead to the two solutions $(p,q,a,b)=(2,5,2,1),(3,7,2,1)$ of (1.2).

Let $p^a=x,q^b=y$ for $(i,j)=(4,2)$ in (1.2), let
\[x=\frac{X+3}{6},~y=Y+2,\] we have
\[Y^2=3X^4+6X^3-3X^2-6X-63.\]

Using \textbf{Magma}, we get all the integral points on this elliptic curve, which are
\[(X,\pm Y)=(-3,2;3),\]
it's easy to see that there are no solutions for $(i,j)=(4,2).$ \hfill $\Box$\\

{\it \textbf{Proof of Theorem 1.3.}}  When $k=2,$ for convenience,
put $p^a=y,q^b=x$ for $(i,j)=(2,3),k=2$ and $p^a=x,q^b=y$ for
$(i,j)=(3,2),k=2$ in (1.2), respectively. Let
\[x=\frac{X+12}{12},~y=\frac{Y+36}{72},\] the converse
transformation is
\[X=12x-12,~Y=36(2y-1).\]
Then, we have
\[Y^2=X^3-144X+22032,~Y^2=X^3-144X-19440,\]respectively.

Using \textbf{Magma}, we get all the integral points on the above
two elliptic curves. For $Y^2=X^3-144X+22032$, the points
$(X,Y)=(36,252),(24,180)$ lead to the solutions
$(p,q,a,b)=(2,2,2,2),(3,3,1,1)$ of (1.2). And $Y^2=X^3-144X-19440$
has no integral points, hence $(1.2)$ has no solutions for
$(i,j)=(3,2)$.

Let $p^a=y,q^b=x$ for $(i,j)=(2,4),k=2$ in (1.2), let
\[x=X,~y=\frac{Y+3}{6},\] we have
\[Y^2=3X^4-18X^3+33X^2-18X+153.\]

Using \textbf{Magma}, we get all the integral points on this elliptic curve, which are
\[(X,\pm Y)=(-1,4;15),\]
by some calculations, we find that $(X,Y)=(4,15)$ lead to the unique
solutions $(p,q,a,b)=(3,2,1,2)$ of (1.2).

Let $p^a=x,q^b=y$ for $(i,j)=(4,2),k=2$ in (1.2), let
\[x=X,~y=\frac{Y+3}{6},\] we have
\[Y^2=3X^4-18X^3+33X^2-18X-135.\]

Using \textbf{Magma}, we get all the integral points on this elliptic curve, which are
\[(X,\pm Y)=(-2,5;15),\]
by some calculations, we find that $(X,Y)=(5,15)$ lead to the unique solution $(p,q,a,b)=(5,3,1,1)$ of (1.2). \hfill $\Box$\\

In figure 1, we display the graph of the quartic curve
\[Y^2=3X^4-18X^3+33X^2-18X-135.\]

\begin{figure}[h!]
\includegraphics[scale=.5]{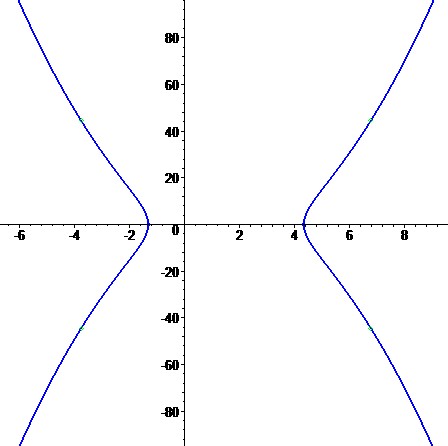}
\caption{$Y^2=3X^4-18X^3+33X^2-18X-135$} \label{hreference labeli}
\end{figure}

\vskip20pt
\section{Two conjectures related with Catalan equation}
In \cite{Cai-Chen}, the first author raised a new variant of the
Hilbert-Waring problem: to express a positive integer $n$ as the sum
of $s$ positive integers whose product is a $k$-th power, i.e.,
\[n=x_1+x_2+\cdots+x_s\]
such that \[x_1x_2\cdots x_s=x^k,\]
for $n,x_i,x,k\in \mathbb{Z}^+$, which may be regarded as a generalization of Waring's problem:
\[n=x_1^k+x_2^k+\cdots+x_s^k.\]

Now we expand this idea to Catalan's equation. Let's consider
\begin{equation}
 \begin{cases}
A-B=1, \\
AB~is~square$-$full,
\end{cases}
\end{equation}
where $A,B$ are positive integers.

This is a generalization of Catalan equation. By using the method of
Pell's equation, it's easy to show that there are infinitely many
solutions of $(4.1)$, the least three are
$(8,9),(288,289),(675,676).$

However, after calculations with computer, we find and raise the
following conjectures (we have checked up to $B<A<10^6$)

\begin{conjecture}\label{conj41}
Let $r\geq0$ be integer, the Diophantine equation
\[
 \begin{cases}
A-B=2^r, \\
AB~is~cube$-$full,
\end{cases} \]
has no solutions for $r=0,$ and has unique solution
$(A,B)=(2^{r+1},2^r)$ for $r\geq1.$
\end{conjecture}

Moreover, we have

\begin{conjecture}\label{conj42}
Let $r\geq1$ be integer, there are infinitely many prime $p$ such
that the Diophantine equation
\begin{equation}
 \begin{cases}
A-B=p, \\
AB~is~cube$-$full,
\end{cases}
\end{equation}
has no solutions. Moreover, the least prime is 29.
\end{conjecture}

It's easy to verify that for every integer $2\leq n\leq 28$, (4.2)
has solutions. However, we even don't know if there is a solution
for infinitely many prime $p$.

\vskip20pt
\section{ Acknowledgements}

The authors wish to thank Dr. Deyi Chen and Mr. Tanyue Gao for their
kind help in calculation, especially in verifying Conjectures 2.1,
4.1 and 4.2.

\vskip20pt
\bibliographystyle{amsplain}

\end{document}